\documentclass[11pt,reqno,twoside]{article}
\usepackage{mathrsfs}
\usepackage{amsthm}
\usepackage{enumitem}
\usepackage[utf8]{inputenc}
\usepackage{fontenc}
\usepackage{hyperref}
\hypersetup{  
   bookmarks=true,
   backref=true,
   pagebackref=false,
   colorlinks=true,
   linkcolor=blue,
   citecolor=red,
   urlcolor=blue
}

\usepackage{subfig}
\usepackage[english]{babel}
\usepackage{amsbsy,amscd,fancyhdr,graphicx,psfrag,fancybox,indentfirst,color}
\usepackage{graphics,epsfig}
\usepackage{amsmath,amsfonts,amssymb}
\usepackage{mysty}
\usepackage{fullpage}
\usepackage{graphicx}

\newcommand*\samethanks[1][\value{footnote}]{\footnotemark[#1]}

\title{\textbf {Continuum limit of the nonlocal $p$-Laplacian evolution problem on random inhomogeneous graphs}}         
\date{}
\author{Hafiene Yosra\thanks{Normandie Univ, ENSICAEN, UNICAEN, CNRS, GREYC, France.} \and Jalal M. Fadili\samethanks \and Christophe Chesneau\thanks{Normandie Univ, ENSICAEN, UNICAEN, CNRS, LMNO, France.} \and Abderrahim Elmoataz\samethanks[1]
}

\begin{document}

\maketitle
\begin{abstract}
In this paper we study numerical approximations of the evolution problem for the nonlocal $p$-Laplacian operator with homogeneous Neumann boundary conditions on inhomogeneous random convergent graph sequences. More precisely, for networks on convergent inhomogeneous random graph sequences (generated first by deterministic and then random node sequences), we establish their continuum limits and provide rate of convergence of solutions for the discrete models to their continuum counterparts as the number of vertices grows. Our bounds reveals the role of the different parameters, and in particular that of $p$ and the geometry/regularity of the data.
\end{abstract}



\section{Introduction}
\subsection{Problem statement}
Our main goal in this paper is to study numerical approximations on random inhomogeneous graphs to a nonlocal nonlinear diffusion problem, involving the nonlocal $p$-Laplacian operator with homogeneous Neumann boundary conditions. More precisely, the nonlocal $p$-Laplacian evolution problem with Neumann boundary conditions that we deal with is
 \begin{equation}\tag{$\mathcal{P}$}
  \begin{cases}
\frac{\partial }{\partial t} u(x,t) =- \A(u(x,t)), \quad x \in \O, t>0,
  \\
  u(x,0)= g(x), \quad x \in \O,
  \end{cases}
  \label{neumann}
  \end{equation}
  where 
\[
\A(u(x,t)) = - \displaystyle{\int_{\O}} K(x,y) \abs{u(y,t) - u(x,t) }^{p-2} (u(y,t) - u(x,t) ) dy ,
\]
with $\O \subset \R $ a compact domain, and without loss of generality $\O = [0,1]$\footnote{Only boundedness of $\O$ is actually needed but we take $\O$ as a closed set as well to conform to our setting of graphs. Moreover, though we here focus on the one-dimensional case $\O \subset \R $, several of our results can be extended to higher dimension.}. The kernel $K \in L^\infty(\O^2)$ is a symmetric and nonnegative mapping. Throughout the paper, we will assume that $p \in ]1,+\infty[$. Existence and uniqueness of a strong solution to \eqref{neumann} in the space $L^p(\O)$ was shown in~\cite[Theorem~3.1]{nonlocal} (relying on arguments from~\cite{Rossimathsurvey2010}).

The interest for this operator has constantly increased over the last few years, as it appears naturally in the study of nonlocal diffusion processes. It arises in a number of applications such as continuum mechanics, phase transition phenomena, population dynamics, image processing and game theory (see \cite{ANDREU2008201,Rossimathsurvey2010,Gilboa2009, kindermann2005} and the references therein). 
On the other hand, recently, there has been a high interest in adapting and applying disecretized versions of PDEs such as \eqref{neumann} on data defined on arbitrary graphs and networks. Given the discrete nature of data in practice, graphs constitute a natural structure suited to their representation. The demand for such methods is motivated by existing and potential future applications, such as in machine learning and mathematical image processing (see among other references \cite{elmoataz2012, elmoataz2015, Gilboa2007, Buades2006}). Indeed, any kind of data can be represented by a graph in an abstract form in which the vertices are associated to the data and the edges correspond to relationships within the data. These practical considerations naturally lead to a discrete time and space approximation of \eqref{neumann}. 
 
To do this, fix $n \in \N^*$. Let  $G_{n} = \pa{V(G_{n}) , E(G_{n})}$, where $V(G_n)$ stands for the set of nodes and $E(G_n) \subset V(G_n) \times V(G_n)$ denotes the edges set, be a sequence of simple graphs, i.e.~undirected graphs without loops and parallel edges.

Next, we  consider the fully discrete counterpart of~\eqref{neumann} on a graph $G_n$ using the forward Euler scheme. For that, let us consider a partition (not necessarily uniform) $\{\tauh\}_{h=1}^{N}$, $N \in \N^*$ of the time interval $[0,T]$ of maximal size $\tau = \max\limits_{h \in  \{1, \cdots, N\}} \tauh $, i.e; $\tauhm \eqdef \aabs{t_h -  t_{h-1}}$. Denote $u_i^h \eqdef u (x_i,t_h)$ and $g_i \eqdef g(x_i)$. Then for $h \in  \acc{1, \cdots, N}$, consider 
\begin{equation}\tag{$\mathcal{P}^{d}_{n, \tau}$}
\begin{cases}
\displaystyle{\frac{\hi - \Hi}{\tauhm}}= \frac{1}{n}\sum\limits_{j: (i,j) \in E(G_n)}\abs{u_{j}^{h-1}- \Hi}^{p-2} (u_{j}^{h-1}- \Hi),
\\ 
u_{i}^0 = g_i , i \in \acc{1, \cdots, n}.
\end{cases}
\label{neumanntotaldiscrete}
\end{equation} 
Thus,~\eqref{neumanntotaldiscrete} induces a discrete diffusion process parametrized by the structure of the graph whise adjacency matrix captures the (nonlocal) interactions. 
As such, it can be viewed as a discrete approximation of a continuous problem such as~\eqref{neumann}.

Several questions then naturally arise: 
\begin{itemize}
\item Does the discrete problem~\eqref{neumanntotaldiscrete}, and in what sense, has a continuum limit (as $n \to +\infty$) ?
\item What is the rate of convergence to this limit ? Is this limit consistent/related with the unique strong
solution of~\eqref{neumann} ?
\item What are the parameters involved in this rate and what is their influence on the convergence rate ?
\end{itemize}

This paper provides answers to these questions for graphs drawn from a random model. The 'classical' random graph models, in particular dense graphs, are 'homogeneous', in the sense that the nodes degrees tend to be concentrated around a typical value, so that  all vertices are exactly equivalent in the definition of the model. Furthermore, in a typical realization, most vertices are in some sense similar to most others. In contrast, many graphs arising in the real world applications do not have this property and are inhomogeneous. One reason is that the vertices may have been 'born' at different times, with old and new vertices having very different properties. In particular, in many examples the degree distribution follows a power law. Thus, there has been a lot of recent interest in defining and studying networks in 'inhomogeneous' random graph models (see Section~\ref{sec:randomgraphmodelrandom} for further details). That is why our aim is to investigate this graph model to study the limit $p$-Laplacian discrete approximation.

\subsection{Contributions and relation to prior work}
In~\cite{Medvedevrandom} and earlier~\cite{medv}, the author studied convergence of discrete approximations of a nonlinear heat equation governed by a Lipschitz continuous potentiel, first on deterministic graphs and then on random ones, both being dense, without discretization of time. This last result can not be applied to the $p$-Laplacian, which requires much more sophisticated arguments. Moreover, the result in~\cite{medv} are asymptotic by nature as they essentially reply on the central limit theorem.
 
In~\cite{nonlocal}, we provided a rigorous justification of the continuum limit~\eqref{neumann} for the discrete $p$-Laplacian on deterministic dense graphs. The analysis of the continuum limit in~\cite{nonlocal} uses ideas from the theory of dense graph limits~\cite{Lovaszdense2006, Lovasz2008, Lovasz2012}, which for every convergent family of dense graphs defines the limiting object, a measurable symmetric and bounded function $K$. This function is called a \textit{graphon}. It captures the connectivity of $G_n$ for large $n$. In~\cite{nonlocal}, for convergent sequences of deterministic dense graphs $\{G_n\}_{n \in \N}$, it was shown that with the kernel in~\eqref{neumann} taken to be the graphon associated to $\{G_n\}_{n\in \N}$, the solution of~\eqref{neumann} is well-approximated by those of the totally discrete problems~\eqref{neumanntotaldiscrete} for large $n$ and small discretization time step $\tau$. However, the analysis in~\cite{nonlocal} does not cover networks on inhomogeneous graphs nor does it deal with random graph models. The latter have many important applications. The main contribution of our paper is to bridge this gap by focusing on evolution systems on inhomogeneous random graphs. 
  
Combining tools from evolution equations, random graph theory and deviation inequalities, we establish \textit{nonasymptotic} rate of convergence of the discrete solution to its continuum limit with high probability. More precisely, we start by considering the case of random graph models generated by a deterministic sequence of nodes. We prove nonasymptotic error bounds that hold with high probability. These results serve as a basis to deal with the totally random graph model, i.e.; where both the nodes and edges are random. In turn, this shows convergence of solutions for the discrete model to the solution of the continuum problem as the number of vertices $n$ grows. To get the corresponding convergence rate, we additionally assume that the kernel $K$ and the initial data $g$ belong to the very large class of to the Lipschitz spaces $\Lip(s, L^q(\O^2))$ and $\Lip(s', L^q(\O))$. Roughly speaking, $\Lip(s, L^q(\O^2))$ contains functions with $s$ "derivatives" in $L^q(\O^2)$. They contain in particular functions of bounded variation and those of fractal structure for appropriate values of $s$, see (see~Appendix~\ref{sec:appendix} for a brief introduction to these functional spaces). Using in addition arguments from approximation theory on these spaces, we get  convergence rates that reveal the role of the value of $p$ and the regularity of the graphon $K$ and the initial data $g$ both on the rate and the probability of success. In particular, we isolate three different regimes where the rate exhibits different scalings.

\subsection{Paper organization}
The rest of the paper is organized as follows. In Section~\ref{sec:randomgraphmodelrandom}, we give the definition of the inhomogeneous random model that we deal with throughout the paper and specify the assumptions needed to get our results. We finish the section by giving an example for which our assumptions are verified. Section~\ref{sec: deterministic} is devoted to the main result of the paper. We begin our analysis by treating random graph sequences generated by deterministic nodes in Section~\ref{subsec: deterministic}. Then, in Section~\ref{subsec : random} we consider the general model defined previously in~Section~\ref{sec:randomgraphmodelrandom}. After getting the convergence of the discrete model to its continuum limit and identifying the corresponding rate, in Section~\ref{subsec-regimes}, we discuss the different regimes of the convergence rate as a function of the problem parameters. Some technical material is deferred to Appendix~\ref{sec:appendix}.
\subsection{Notations}
For a graph $G = \pa{V(G) , E(G)}$, two vertices $i, j \in V(G)$ are adjacent, if they are connected by an edge. Let $G_n = \pa{V(G_n),E(G_n)}$, $n \in  \N^*$, be a sequence of inhomogeneous, finite, and simple graphs.
 
 For a given vector $u = (u_1,\cdots, u_n )^{\top} \in \R^n$, we define the norm $\anorm{\cdot}_{p,n}$ 
 \[
 \anorm{u}_{p,n}=  \left ( \frac{1}{n} \sum\limits_{i=1}^{n} \abs{u_i}^p\right )^{\frac{1}{p}}.
 \]
For an integer $n \in \N^*$, we denote $[n]=\{1,\cdots,n\}$. For any set $S$, $\cl{S}$ is its closure and $\aabs{S}$ is its cardinality or its Lebesgue measure (to be understood from the context). $\chi_{S}$ is the characteristic function of the set $S$ (takes $1$ in it and $0$ otherwise).

$C(0,T; L^p(\O))$ denotes the space of uniformly time continuous functions with values in $L^p(\O)$. For $d \in \acc{1,2}$, $\Lip(s,L^q(\O^d))$ is the Lipschitz space which consists of functions with, roughly speaking, $s$ "derivatives" in $L^q(\O^d)$~\cite[Ch.~2, Section~9]{devorelorentz93}. Only values $s \in ]0,1]$ are of interest to us. See Section~\ref{subsec:spaces} for further details on these spaces and approximation theoretic results on them.

\section{The random inhomogeneous graph model}
\label{sec:randomgraphmodelrandom}

\subsection{The graph model}
We start with the description of the model of inhomogeneous random graphs that will be used throughout. This random
graph model is motivated by the construction of inhomogeneous random graphs in~\cite{Ballobas2006,Ballobas2008,Ballobaslong2008}.

\begin{defi}
Fix $n \in \N^*$ and let $K$ be a symmetric measurable function on $\O^2$. Generate the graph $G_n = \pa{V(G_n),E(G_n)} \eqdef G_{q_n}(n,K)$ as follows: 
\begin{enumerate}[label=\arabic*)]
\item Generate $n$ independent and identically distributed (i.i.d.) random variables $(\bX_1, \cdots, \bX_n) \eqdef \bX$ from the uniform distribution on $\O$. Let $\acc{\bX_{(i)}}_{i = 1}^{n}$ be the order statistics of the random vector $\bX$, i.e. $\bX_{(i)}$ is the $i$-th smallest value.
\item Conditionally on $\bX$, join each pair $(i, j) \in [n]^2$ of vertices independently, with probability $q_n \wedgX$, i.e. for every $(i,j) \in [n]^2$, $i \neq j$,
\begin{equation}
 \P\pa{(i,j) \in E(G_n) | \bX} = q_n \wedgX , 
 \label{def:graphmodel}
\end{equation}
where 
\begin{equation}
\wedgX \eqdef \min \pa{\frac{1}{\abs{\OXij}}\int_{\OXij}K(x, y) dx dy, 1/q_n},
\label{def:wedgX}
\end{equation}
and 
\begin{equation}
\OXij \eqdef ]\bX_{(i-1)}, \bX_{(i)}]  \times ]\bX_{(j-1)}, \bX_{(j)}] 
\end{equation}
where $q_n$ is non-negative and uniformly bounded in $n$. 
\end{enumerate}
A graph $G_{q_n}(n,K)$ generated according to this procedure is called a $K$-random inhomogeneous graph generated by a random sequence $\bX$.
\label{def : randomgraph}
\end{defi}

At this stage, the following important remark is in order.
\begin{rem}
In the context of numerical analysis, we are primarily interested not only in the error bounds of the discrete problem, but more importantly in the (nonasymptotic) rate of convergence. This is why our attention aims specifically at this graph model and not at the original inhomogeneous random model defined in~\cite{Ballobas2006, Ballobas2008}, i.e.~the model constructed replacing~\eqref{def:graphmodel} by 
\begin{equation*}\label{eq:originhom}
\P\pa{(i,j) \in E(G_n)} = \min\pa{q_nK(\bX_{i}, \bX_{j}), 1}.
\end{equation*}
Our error bounds of the discrete problem~\eqref{neumanntotaldiscrete} cover also this graph model, and more specifically, the first statements of Theorem~\ref{theo : convergencedeterministe} and Theorem~\ref{theo : totalrandom} hold. However, with this model, even our convergence claim (not to mention the rate) of the discrete scheme does not hold unless the kernel $K$ and the intial data $g$ are additionally supposed almost everywhere continuous.
\label{rem:classicalmodel}
\end{rem}

We denote by $\bx = (\bx_1, \cdots, \bx_n)$  the realization of $\bX$. To lighten the notation, we  also denote 
\begin{equation}
\O^{\bX}_{ni} \eqdef ]\bX_{(i-1)}, \bX_{(i)}], \quad \Oxi \eqdef   ]\bx_{(i-1)}, \bx_{(i)}], \quad \text{and} \quad   \Oxij \eqdef  ]\bx_{(i-1)}, \bx_{(i)}] \times   ]\bx_{(j-1)}, \bx_{(j)}] \quad  i,j \in [n].
\label{intervalpartition}
\end{equation}
As the realization of the random vector $\bX$ is fixed, we define 
\begin{equation}
\wedg \eqdef \min\pa{\frac{1}{\abs{\Oxij}}\displaystyle{\int_{\Oxij}}K(x, y) dx dy, 1/q_n}, \quad \forall (i,j) \in [n]^2, \quad  i \neq j.
\label{weightedgraph}
\end{equation}
In the rest of the paper, the following random variables will be useful. Let $\lam$, $(i,j) \in [n]^2, i\neq j$, be i.i.d. random variables such that $q_n \lam$ follows a Bernoulli distribution with parameter $q_n \wedg$. We consider the i.i.d. random variables $\ups$ such that the distribution of $q_n \ups$ conditionally on $\bX=\bx$ is that of $q_n \lam$. Thus $q_n \ups$ follows a Bernoulli distribution with parameter $\EE\bpa{q_n \wedgX}$, where $\EE(\cdot)$ is the expectation operator (here with respect to the distribution of $\bX$). 

\vskip 12pt

We now formulate our assumptions on the graph sequence $\acc{G_{q_n}(n,K)}_{n \in \N}$.
\begin{assum} We suppose that $q_n$ and $K$ are such that the following hold: 
\begin{enumerate}[label=\rm ({A}.\arabic*)]
\item $G_{q_n}(n,K)$ converges almost surely and its limit is the graphon $K \in L^{\infty} (\O^2)$;\label{assum:A1}
\item $\inf\limits_{n \geq 1} q_n > 0$ and $\sup\limits_{n \geq 1} q_n < + \infty$.\label{assum:A2}
\end{enumerate}
\label{assumption}
\end{assum}

\subsection{Example} 
\label{subsec:graphexample}
Although we shall give a general result throughout the paper, it may help to bear in mind one particular example of the general class of models we shall study. This example is inspired by the so-called \textit{almost dense (or non uniform)} random graphs (see \cite[Section~3.4]{Ballobas2008}).

\begin{prop} Suppose $K  \in L^{\infty} (\O^2)$ is a symmetric measurable function. Choose the parameter $q_n = n^{- g(n)}$ such that $g(n) \log(n) = O(1)$. Then, assumptions \ref{assum:A1} and \ref{assum:A2} are in force. 
\end{prop}
\bpf{}  
Since the graphon $K \in L^{\infty}(\O^2)$, the arguments to prove~\cite[Lemma~3.5 and Lemma~3.8]{Ballobas2008}, that were designed for the graph model described in Remark~\ref{rem:classicalmodel}, can be adapted to cover our graph model with~\eqref{def:graphmodel} to show that the sequence of random graphs $G_{q_n}(n, K)$ indeed converges almost surely to the graphon $K$ in the metric $d_{\text{sup}}$ (see \cite[Section~2.1]{Ballobas2008} for details about this metric). This shows~\ref{assum:A1}.
As we suppose that $g(n) \log (n) = O(1) $, we get immediately that~\ref{assum:A2} is verified.
\epf

Observe that taking the trivial choice $q_n = O(1)$, one recovers the dense random graph model extensively studied in~\cite{Lovaszdense2006, Lovaszrandom2011}.

\section{Consistency of the nonlocal $p$-Laplacian on random inhomogeneous graphs}
\label{sec: deterministic}
Having defined the structure of the network, we are now in position to state our main error bounds between the discrete dynamics and their continuous ones. First, in Section~\ref{subsec: deterministic}, we assume that $\bX$ is deterministic. Capitalizing on this result, we will then deal with the totally random model (i.e.; generated by random nodes) in Section~\ref{subsec : random} by a simple marginalization argument.

\subsection{Networks on graphs generated by deterministic nodes}
\label{subsec: deterministic}
 
We define the parameter $\delta(n)$ as the maximal size of the spacings between the the ordered values $\bx_{(i)}$
\begin{equation}
\delta(n) = \max\limits_{i \in [n]} \abs{\bx_{(i)} - \bx_{(i-1)}}.
\label{spacingdeterministic}
\end{equation}

Next, we consider the following system of difference equations on $G_{q_n}(n,K)$\footnote{This is clear by proper normalization by $q_n$ (by dividing and multiplying by $q_n$). We abuse notation to lighten the system.} : 
\begin{equation}\tag{$\mathcal{P}^{d,d}_n$}
\begin{cases}
\displaystyle{\frac{\hi - \Hi}{\tauhm}}= \frac{1}{n}\sum\limits_{j = 1}^{n} \lam\abs{u_{j}^{h-1}- \Hi}^{p-2} (u_{j}^{h-1}- \Hi), \quad (i,h) \in [n] \times [N],
\\ 
u_{i}^0 = g_i , i \in [n],
\end{cases}
\label{neumanndeterministicdiscrete}
\end{equation} 
where
\[
g_i = \frac{1}{\abs{\Oxi}}\int_{\Oxi} g(x) dx .
\]
Recall from Section~\ref{sec:randomgraphmodelrandom} that $\lam$ are the i.i.d. random variables such that $q_n \lam $ follows the Bernoulli distribution with parameter $q_n \wedg$.

Before turning to our convergence result, we pause here to make the following two important observations. 
\begin{rem} 
Coming back to Definition~\ref{def : randomgraph}, one can easily check that $G_{q_n}(n, K)$ is actually a product probability space\footnote{To keep notation simple, we allow for loops, in our random graph model. Excluding loops would not lead to any changes in the analysis.}
\[
\boldsymbol{\O}_n \eqdef \boldsymbol{\O}_n^V \times \boldsymbol{\O}_n^E  \eqdef  \pa{\O_n^V \eqdef [0,1]^{n}, 2^{\O_n^V}, \mathbb{P}} \times \pa{\O_n^E \eqdef \acc{0,1}^{n (n + 1)/2}, 2^{\O_n^E}, \mathbb{P}} .
\]
So that, rigorously speaking, if we take a random event $\boldsymbol{\omega}$ from $\boldsymbol{\O}_n$, problem~\eqref{neumanndeterministicdiscrete} must be written using $\lam (\boldsymbol{\omega})$ instead of $\lam$, and likewise for all other random variables. For notational simplicity, we drop $\boldsymbol{\omega}$. But it is important to keep in mind that the evolution equations we write involving random variables must be understood in this sense. 
\label{rem:event}
\end{rem}

\begin{rem} 
As the reader may have remarked, the sum in the right-hand side of~\eqref{neumanndeterministicdiscrete} is divided by $n$ instead of a weighted sum with weights $\aabs{\bx_{(i)} - \bx_{(i-1)}}^{-1}$ which would be expected if we interpret this sum as a Riemann sum. The scaling by $n$ reminds us of an equidistant design regarding the space-discretization, despite the fact that the nodes are chosen not necessarily equispaced. However, given that the $\bx_i$'s are realizations of i.i.d. uniform variables on $\O$, the uniform spacing choice still makes sense. Indeed, using classical results on order statistics of uniform variables, see, e.g.,~\cite[Section~1.7]{Reis1989}, it can be shown that each spacing $\bX_{(i)} - \bX_{(i-1)}$ concentrates around $i/n$ for $i \in [n]$.
\end{rem} 

We are now in position to tackle our main goal: comparing the solutions of the discrete and continuous problems and establish our rate of convergence. Since the two solutions do not live on the same spaces, it is reasonable to represent some intermediate model that is the continuous extension of the discrete problem, using the vector $U_h = (u_1^h, u_2^h, \cdots , u_n^h)^{\top}$ whose components uniquely\footnote{In \cite[Lemma 5.1]{nonlocal}, we show that \eqref{neumanndeterministicdiscrete} is well posed.} solve the previous system~\eqref{neumanndeterministicdiscrete} to obtain the following piecewise linear interpolation on $\O \times  [0, T ]$
\begin{equation}
\U(x,t) = \frac{t_h - t }{\tauhm} \Hi + \frac{t-t_{h-1}}{\tauhm} \hi \quad  \text{if} \quad x \in \Oxi, \quad t\in ]t_{h-1},t_h],
\label{eq:extensiondet}
\end{equation}
and a piecewise approximation 
\begin{equation}
\M (x,t) = \sum_{i=1}^{n}\sum_{h=1}^N \Hi \chi_{]t_{h-1},t_h]} (t) \chi_{\Oxi}(x) .
\label{eq:piecewisedet}
\end{equation}
Then, $\U$ uniquely solves the following problem
\begin{equation}\tag{$\mathcal{P}_n$}
\begin{cases}
\frac{\partial}{\partial t} \U(x,t)=-\boldsymbol{\Delta}^{\Lam}_p(\bar{u}_n(x,t)), \quad x \in \O, t>0,\\
\U(x,0) = g_n(x),\quad x \in \O ,
\end{cases}
\label{neumannenxtension}
\end{equation}
where the random variable 
\[
\Lam(x,y) = \lam \quad \text{for} \quad (x,y) \in \Oxij,
\]
and 
\[
g_n(x) = g_i \quad \text{if} \quad x \in \Oxi, i \in [n].
\]
Toward our goal of establishing error bounds, we need an intermediate discrete problem for the $p$-Laplacian. This is defined as
\begin{align}\tag{$\stackrel{\wedge}{\mathcal{P}}_n^{d}$}
\begin{cases}
\displaystyle{\frac{v_i^h - v_i^{h-1}}{\tauhm}}= \frac{1}{n}\sum\limits_{j = 1}^{n} \wedg \abs{v_{j}^{h-1}- v_{i}^{h-1}}^{p-2} (v_{j}^{h-1}- v_{i}^{h-1}), \quad (i,h) \in [n] \times[N], \\ 
v_{i}^0 = g_i , \quad i \in [n].
\end{cases}
\label{neumanndiscrete2}
\end{align}
The discrete problem~\eqref{neumanndiscrete2} can also be viewed as a discrete $p$-Laplacian evolution problem over a complete\footnote{Recall that a complete graph is a simple undirected graph in which each pair of vertices is connected by an edge.} weighted graph on $n$ vertices, where the weight of edge $(i,j)$ is $\wedg$.

Using the vector~$V_n^h = (v_1^h,v_{2}^h,\cdots, v_{n}^h)^{\top}$ whose components uniquely solve the system~\eqref{neumanndiscrete2} , similarly to before, we define the following linear interpolation on $\O \times [0,T]$
\begin{equation}
\check{v}_n(x,t) = \frac{t_h - t }{\tauhm}v_i^{h-1}+ \frac{t-t_{h-1}}{\tauhm} v_i^{h} \quad  \text{if} \quad x \in \Oxi,\quad t\in ]t_{h-1},t_h],
\label{eq:extensiondet1}
\end{equation}
and a piecewise-constant approximation
\begin{equation}
\bar{v}_n(x,t) =  \sum_{i=1}^{n}\sum_{h=1}^N v_i^{h-1}  \chi_{]t_{h-1},t_h]} (t) \chi_{\Oxi}(x).
\label{piecewise1}
\end{equation}
We also define the piecewise-constant extension $\Wedg$ on $\O^2$ 
\begin{equation}
\label{eq:contextcheckKn}
\Wedg(x,y) = \sum_{(i,j) \in [n]^2} \wedg \chi_{\Oxij}(x,y) .
\end{equation}
Then, by construction, $\check{v}_n(x,t) $ uniquely solves the following problem 
\begin{equation}\tag{$\stackrel{\wedge}{ \mathcal{P}_n}$}
\begin{cases}
\frac{\partial}{\partial t} \check{v}_n(x,t)=-\boldsymbol{\Delta}^{\Wedg}_p(\bar{v}_n(x,t)), \quad x \in \O, t>0,\\
\check{v}_n(x,0)= g_n(x), \quad x \in \O,
\end{cases}
\label{neumann1}
\end{equation}
where 
\[
g_n(x)=  g_i  \quad \text{for} \quad x \in \Oxi, i \in [n].
\]
 
The first main result of the paper is the following theorem.
\begin{theo}
Suppose that $p \in ]1,+\infty [$, $K \in L^\infty(\O^2)$ is a symmetric and measurable mapping, and $g \in L^{\infty}(\O)$. Let $u$ and $U_h$ denote the unique solutions to~\eqref{neumann} and~\eqref{neumanndeterministicdiscrete}, respectively. Let $\U$ be the continuous extension of $U_h$ given in~\eqref{eq:extensiondet}. Then, the following hold:
\begin{enumerate}[label=(\roman*)]
\item for $T> 0$, there exists a positive constant $C$, such that for any $\beta > 0$
{\small
\begin{equation}
\hspace*{-1cm}
\anorm{u - \U}_{C(0,T; L^p(\O))} \le C \pa{\beta \frac{\log (n)}{n} + \frac{\max\pa{q_n^{-(p-1)}, q_n^{-p/2}}}{ n^{p/2}}}^{1/p} + \anorm{K - \Wedg}_{L^p(\O^2)} +\anorm{g - g_n}_{L^p(\O)}+  O(\tau),
\label{eq:estimate}
\end{equation}}
with probability at least $1 - n^{-C  \min\pa{q_n^{2p-1}, q_n^p}\beta}$.
\item Suppose furthermore that  $g \in  \Lip(s, L^q(\O))$ and $K \in \Lip(s', L^q(\O^2))$, $q \in [1, + \infty]$, $s, s' \in ]0,1]$, and $q_n\norm{K}_{L^\infty(\O^2)} \leq 1$. Then, for $T> 0$, there exists a positive constant $C$, such that for any $\beta > 0$
{\small 
\begin{equation}
\anorm{u - \U}_{C(0,T; L^p(\O))} \le C \pa{\pa{\beta \frac{\log (n)}{n} + \frac{\max\pa{q_n^{-(p-1)}, q_n^{-p/2}}}{ n^{p/2}} }^{1/p} + \delta(n)^{\min\pa{s, s'} \min\pa{1, q/p}} }+O(\tau),
\label{eq:estimate2}
\end{equation}}
with probability at least $1 - n^{-C  \min\pa{q_n^{2p-1}, q_n^p}\beta}$, where $\delta(n)$ is the parameter defined in~\eqref{spacingdeterministic}.
\end{enumerate}
\label{theo : convergencedeterministe}
\end{theo}

Before proceeding to the proof, some remarks are in order.
\begin{rem}
$~$\\\vspace*{-0.5cm}
\begin{enumerate}[label=(\roman*)]
\item By Lemma~\ref{lem:quasidist}, it is clear that the first term in the bounds~\eqref{eq:estimate}-\eqref{eq:estimate2} can be replaced by
\[
\beta^{1/p} \pa{\frac{\log (n)}{n}}^{1/p} + \frac{\max\pa{q_n^{-(1-1/p)}, q_n^{-1/2}}}{ n^{1/2}} .
\]
\item The constant in~\eqref{eq:estimate} depends on $p$ and the data via $\norm{g}_{L^\infty(\O)}$ and $\norm{K}_{L^\infty(\O)}$. For the bound~\eqref{eq:estimate2}, it also depends on $(q,s,s')$.
\item One may wonder if the functional space assumption made on $g$ and $K$ in claim (ii) is reasonable or even makes sense. The answer is affirmative. Indeed, Lipschitz spaces are rich enough to include both functions with discontinuities and even fractal structure. For instance, from~\cite{Lovasz2012}, one can show that the graphon corresponding to the nearest neighbour graphs, which are very popular in practice (e.g.~in image processing~\cite{elmoataz2015,elmoataz2012}), are typical examples satisfying Assumptions~\ref{assum:A1}-\ref{assum:A2} with $q_n=1$ and $K$ is a $\acc{0,1}$-valued function living on the space of bounded variation functions, which in turn is $\Lip(1,L^1(\O^2))$.
\end{enumerate}
\label{rem:Cparam}
\end{rem}

To prove Theorem~\ref{theo : convergencedeterministe}, we first show the following key lemma.
\begin{lem} 
Let 
\[
\varepsilon =  \pa{\beta \frac{\log(n)}{n}+C_3  \max\pa{q_n^{-(p-1)}, q_n^{-p/2}}\frac{1}{n^{p/2}}}^{1/p} + O(\tau) .
\]
Under the assumptions of Theorem~\ref{theo : convergencedeterministe}, for $T> 0$, there exists a positive constant $C$, such that for any $\beta > 0$
\begin{equation}
\P\pa{ \anorm{\check{v}_n - \U}_{C(0,T; L^p(\O))}  \ge \varepsilon } \le n^{- C \min\pa{q_n^{2p-1}, q_n^p}\beta} ,
\label{proba}
\end{equation}
\label{bernstein}
(the constant $C_3$ is given in the proof).
\end{lem}

\bpf{of Lemma~\ref{bernstein}}   
For $1< p<+ \infty$, we define the function 
  \begin{equation*}
  \begin{split}
  &\F :   \mathbb{R} \rightarrow \mathbb{R} \\
 &x \mapsto \aabs{x}^{p-2} x = \sign(x) {\aabs{x}}^{p-1}.
  \end{split}
 \end{equation*}
Observe that $\check{v}_n(\cdot,t)$ and $\U(\cdot,t)$ are both constants over $\Oxi$.  Similarly, $\bar{v}_n(\cdot, t ) $ and $\bar{u}_n (\cdot,t) $ are also constants over the cell $\Oxi$. We therefore used the shorthand notations for the vector-valued functions $\bar{{\mathbf{u}}}_{n}(t)=(\bar{{\mathbf{u}}}_{ni}(t))_{i \in [n]} \eqdef (\bar{u}_n (\bx_i,t))_{i \in [n]}$ and $\bar{{\mathbf{v}}}_{n}(t)=(\bar{{\mathbf{v}}}_{n}(t))_{i \in [n]} \eqdef (\bar{v}_n (\bx_i,t))_{i \in [n]}$, and likewise for $\check{{\mathbf{u}}}_{n}(t)$ and $\check{{\mathbf{v}}}_{n}(t)$. Let us denote $\check{\xi}_n(t) = \check{{\mathbf{u}}}_{n}(t) -  \check{{\mathbf{v}}}_{n}(t)$ and $\bar{\xi}_n(t) = \bar{{\mathbf{u}}}_{n}(t)- \bar{{\mathbf{v}}}_{n}(t)$.
By subtracting both sides of~\eqref{neumannenxtension} from those of \eqref{neumann1}, evaluated at the cell $\Oxi$, we obtain 
\begin{equation}
\frac{d}{dt} \qi(t) = Z_{ni}(t) + \frac{1}{n} \sum\limits_{j=1}^{n} \lam \bpa{\F(\bar{{\mathbf{u}}}_{nj}(t) - \bar{{\mathbf{u}}}_{ni} (t))-\F(\bar{{\mathbf{v}}}_{nj}(t) - \bar{{\mathbf{v}}}_{ni} (t))  },
\label{difference1}
\end{equation}
where
\begin{equation}
Z_{ni}(t) = \frac{1}{n} \sum\limits_{j=1}^{n} (\lam - \wedg) \F(\bar{{\mathbf{v}}}_{nj}(t) - \bar{{\mathbf{v}}}_{ni} (t)).
\label{zed}
\end{equation}
For notational convenience, we denote $\alpha_{ij}(t) \eqdef  \F(\bar{{\mathbf{v}}}_{nj}(t) - \bar{{\mathbf{v}}}_{ni} (t))$, for $(i,j) \in [n]^2$, $t \in [0,T]$. We multiply both sides of~\eqref{difference1} by $\frac{1}{n}\F(\check{{\xi}}_{ni}(t)) $ and sum over $i$ to obtain
\begin{equation}
\frac{1}{p} \frac{d}{dt} \anorm{\check{\xi}_{n}(t)}^p_{p,n} = \frac{1}{n}  \sum\limits_{i=1}^{n} Z_{ni}(t)\F(\check{{\xi}}_{ni}(t)) +  \frac{1}{n^{2}} \sum\limits_{i,j=1}^{n}  \lam  \bpa{\F(\bar{{\mathbf{u}}}_{nj}(t) - \bar{{\mathbf{u}}}_{ni}(t))-\F(\bar{{\mathbf{v}}}_{nj}(t) - \bar{{\mathbf{v}}}_{ni}(t)) }\F(\check{{\xi}}_{ni}(t)) .
\label{differencee1}
\end{equation}
We estimate the first term on the right-hand side of~\eqref{differencee1} using the H\"older inequality, to get
\begin{equation}
\frac{1}{n} \aabs{ \sum\limits_{i=1}^{n} Z_{ni}(t) \F(\check{{\xi}}_{ni}(t)) }\le \frac{1}{n} \pa{  \sum\limits_{i=1}^{n}  \abs{Z_{ni}(t)}^p}^{\frac{1}{p}} \times \pa{  \sum\limits_{i=1}^{n} \abs{\check{\xi}_{ni}(t)}^p}^{\frac{p-1}{p}} \le \anorm{Z_n(t)}_{p,n}  \anorm{\check{\xi}_n(t)}^{p-1}_{p,n}.
\label{first1}
\end{equation}

Now, using the fact that $0 \le\lam \le 1/{q_n}$, $\forall (i,j) \in  [n]^2$ and applying~\cite[Corollary~B.1]{nonlocal} to the function $\F $ between $a=\bar{{\mathbf{v}}}_{nj} (t) - \bar{{\mathbf{v}}}_{ni}(t) $ and $b= \bar{{\mathbf{u}}}_{nj}(t) - \bar{{\mathbf{u}}}_{ni}(t)$ (without loss of generality, we suppose that $b >a $), we get 
\begin{equation}
\begin{aligned}
&\aabs{\frac{1}{n^{2}}  \sum\limits_{i,j=1}^{n}  \lam  \bpa{ \F(\bar{{\mathbf{u}}}_{nj}(t) - \bar{{\mathbf{u}}}_{ni}(t)) - \F(\bar{{\mathbf{v}}}_{nj}(t) - \bar{{\mathbf{v}}}_{ni}(t)} \F({\xi}_{ni}(t))} \\
&\le \frac{1}{q_n} \frac{p-1}{n^{2}} \sum\limits_{i,j=1}^{n}  \abs{\bar{\xi}_{nj} - \bar{\xi}_{ni}} \abs{\eta_n(t)}^{p-2} \abs{\check{\xi}_{ni}}^{p-1},
\label{ineq:mean}
\end{aligned}
\end{equation}
where $\eta_n(t)$ is an intermediate value between $a$ and $b$.  Using that fact that $g \in L^{\infty}(\O)$ and the construction of $\bar{{\mathbf{u}}}_{n}(\cdot)$, we deduce from~\cite[Theorem~3.1(ii)]{nonlocal} that for $t \in [0,T]$ 
\begin{equation}
\abs{\eta_n(t)}^{p-2} \leq  \abs{\bar{{\mathbf{u}}}_{nj}(t) - \bar{{\mathbf{u}}}_{ni}(t)}^{p-2} \leq \pa{ 2 \anorm{u(\cdot,t)}_{L^{\infty}(\O)}}^{p-2} \leq \pa{ 2 \anorm{g}_{L^{\infty}(\O)}}^{p-2} = C_2 . 
\label{ineq:meanvalue}
\end{equation}
Inserting~\eqref{ineq:meanvalue} into~\eqref{ineq:mean}, and then using the H\"older and triangle inequalities, it follows that
\begin{equation}
\begin{aligned}
&\aabs{\frac{1}{n^{2}}  \sum\limits_{i,j=1}^{n}  \lam \bpa{ \F(\bar{{\mathbf{u}}}_{nj}(t) - \bar{{\mathbf{u}}}_{ni}(t))- \F(\bar{{\mathbf{v}}}_{nj}(t) - \bar{{\mathbf{v}}}_{ni}(t) }\F(\check{{\xi}}_{ni}(t)) }\\
& \le \frac{C_2}{q_n} \frac{p-1}{n^2} \sum\limits_{j=1}^{n} \sum\limits_{i=1}^{n} \abs{\bar{\xi}_{nj}(t) - \bar{\xi}_{ni}(t)} \abs{\check{\xi}_{ni}}^{p-1}\\
& \le \frac{C_2}{q_n} \frac{p-1}{n^2}\pa{\pa{ \sum\limits_{i,j}  \abs{\bar{\xi}_{nj}(t) - \bar{\xi}_{ni}(t)}^p}^{\frac{1}{p}} \pa{\sum\limits_{i,j}\abs{\check{\xi}_{ni}(t)}^p}^{\frac{p-1}{p}} } \\
& \le  \frac{C_2}{q_n} \frac{p-1}{n^2} \pa{\pa{ \sum\limits_{i,j}\abs{\bar{\xi}_{nj}(t)}^p}^{\frac{1}{p}} +\pa{ \sum\limits_{i,j}\abs{\bar{\xi}_{ni}(t)}^p}^{\frac{1}{p}}  } \pa{ n^{\frac{2(p-1)}{p}} \pa{ \frac{1}{n} \sum\limits_{i=1}^{n} \abs{\check{\xi}_{ni}(t)}^p}^{\frac{p-1}{p}}} \\
& \le \frac{C_2}{q_n} \frac{p-1}{n^2} \pa{ 2 n^{\frac{2}{p}}  \anorm{\bar{\xi}_n(t)}_{p,n}  }  \pa{ n^{\frac{2(p-1)}{p}} \anorm{\check{\xi}_n(t)}^{p-1}_{p,n}  } \\
& \le2 C_2 \frac{p-1}{q_n} \anorm{\bar{\xi}_n(t)}_{p,n}\anorm{\check{\xi}_n(t)}^{p-1}_{p,n}.
\end{aligned}
\label{second1}
\end{equation}
 Using the triangle inequality combined with the result of~\cite[Lemma~5.2]{nonlocal}, we have
\begin{equation}
\begin{aligned}
\anorm{\bar{\xi}_n(t)}_{p,n} &=   \anorm{\bar{{\mathbf{v}}}_n(t) - \bar{{\mathbf{u}}}_n(t)}_{p,n}\\
& \le \anorm{\bar{{\mathbf{v}}}_n(t) - \check{{\mathbf{v}}}_n (t)}_{p,n} + \anorm{\check{{\mathbf{v}}}_n(t) - \check{{\mathbf{u}}}_n (t)}_{p,n}  + \anorm{\check{{\mathbf{u}}}_n(t) - \bar{{\mathbf{u}}}_n (t)}_{p,n} \\
&\le C \tau + \anorm{\check{\xi}_n(t) }_{p,n}  + C' \tau\\
&\le C'' \tau + \anorm{\check{\xi}_n(t) }_{p,n}.
\end{aligned}
\label{intermidiaire}
\end{equation} 
 
Putting together~\eqref{first1},~\eqref{second1} and \eqref{intermidiaire}, we have 
\begin{equation}
\begin{aligned}
\frac{d}{dt} \anorm{\check{\xi}_n(t)}^p_{p,n} &\le  \anorm{Z_n(t)}_{p,n} \anorm{\check{\xi}_n (t)}_{p,n}^{p-1}+2 C_2 (p-1)/ q_n\pa{C'' \tau + \anorm{\check{\xi}_n(t) }_{p,n} }  \anorm{\check{\xi}_n (t)}_{p,n}^{p-1}\\
&\le \pa{  2 C_3 (p-1) \tau +\anorm{Z_n(t)}_{p,n} } \anorm{\check{\xi}_n(t) }_{p,n}^{p-1} +  2 C_2 (p-1) /q_n \anorm{\check{\xi}_n(t) }_{p,n}^{p}.
\end{aligned}
\label{ineq}
\end{equation}
Then, from~\eqref{ineq} via the Gronwall's inequality in its differential form (see, e.g.,~\cite[Appendix~B]{Evans2010}), we obtain 
\begin{equation}
\sup_{t \in [0,T]} \anorm{\check{\xi}_n(t)}_{p,n} \le \pa{ 2 C_3 T \tau + \sup_{t \in [0,T]} \anorm{Z_n(t)}_{p,n} } \exp \left \{\frac{2 C_2 T}{q_n} \right \}.
\label{ineq:gronwall}
\end{equation}
Since we suppose that $q_n$ verifies Assumption~\ref{assum:A2}, then $\exp \left \{\frac{2 C_2 T}{q_n} \right \}$  is a bounded quantity.
It remains to bound $\sup\limits_{t \in [0,T]} \anorm{Z_n(t)}_{p,n}$. For this purpose, we use Lemma~\ref{lem:process} (see Section~\ref{subsec:devZ})\footnote{This inequality is sharp as can be seen for instance from assertion~(ii) of Lemma~\ref{lem:process}, at leat for $p \geq 2$.}.
Thus, plugging~\eqref{ineq:probaZ} into inequality~\eqref{ineq:gronwall}, we get the desired conclusion. 
\epf

We are now ready to prove our main result.
\bpf{of Theorem~\ref{theo : convergencedeterministe}} 
\begin{enumerate}[label=(\roman*)]
\item   Using the triangle inequality, we have 
  \begin{equation}
  \anorm{u - \U}_{C(0,T; L^p(\O))} \le  \anorm{u - \check{v}_n}_{C(0,T; L^p(\O))} +  \anorm{ \check{v}_n - \U}_{C(0,T; L^p(\O))}.
  \label{ineq:triangle}
  \end{equation}

Since by construction $\Wedg$ is a bounded mapping, we bound the first term on the right-hand side of~\eqref{ineq:triangle} using~\cite[Theorem~5.1]{nonlocal} to get
\begin{equation}
\anorm{u - \check{v}_n}_{C(0,T; L^p(\O))} \leq \anorm{K - \Wedg}_{ L^p(\O^2)} +  \anorm{g - g_n}_{ L^p(\O)}+O(\tau).
\label{ineq: estimationwedge}
\end{equation}
Inequality~\eqref{eq:estimate} then follows by combining~\eqref{ineq: estimationwedge} with~\eqref{proba}.
\item Our assumption on $q_n$ together with~\eqref{weightedgraph} and~\eqref{eq:contextcheckKn} entail that 
\[
\Wedg(x,y) = \sum_{(i,j) \in [n]^2} K_{nij} \chi_{\Oxij}(x,y), \quad K_{nij} = \frac{1}{\abs{\OXij}}\int_{\OXij} K(x,y) dx dy
\] 
Since $g \in \Lip (s, L^q(\O))$ and $K \in \Lip (s', L^q(\O^2))$, we can invoke Lemma~\ref{lem:estimateglip} to get 
\begin{equation}
\anorm{K- \Wedg}_{ L^p(\O^2)} \leq C(p,q,s') \delta(n)^{s' \min\pa{1, q/p}} \qandq \anorm{g - g_n}_{L^p(\O)} \leq C(p,q,s) \delta(n)^{s \min\pa{1, q/p}} .
\label{eq: LipK}
\end{equation}
Inserting the bound~\eqref{eq: LipK} into~\eqref{eq:estimate}, and using the fact that $\delta(n) < 1$, yields~\eqref{eq:estimate2}.
\end{enumerate}
\epf

\subsection{Networks on graphs generated by random nodes}
\label{subsec : random}
Let us now turn to the totally random graph model. 
Consider the following system of difference equations on the totally random graph $G_{q_n}(n, K)$\footnote{\label{foot:drop}Recall again from Remark~\ref{rem:event}, that rigorously speaking, each random variable involved in the problems and equations of this section should be understood as a function of an event $\boldsymbol{\omega}$ from $\boldsymbol{\O}_n$. This dependence is dropped only to lighten notation.} :
\begin{equation}\tag{$\mathcal{P}^{r,d}_n$}
\begin{cases}
\displaystyle{\frac{\hi - \Hi}{\tauhm}}= \frac{1}{n}\sum\limits_{\ens{j}{(i,j) \in E(G_{q_n}(n,K))}}\abs{u_{j}^{h-1}- \Hi}^{p-2} (u_{j}^{h-1}- \Hi), \quad h \in [N]
\\ 
u_{i}^0 = g_i , i \in [n].
\end{cases}
\label{neumannrandomdiscrete}
\end{equation} 
As we have done before, we consider the continuous extension of the solution vector $U^h = (u_1^h, \linebreak u_2^h, \cdots, u_n^h)^\top$, that is a linear interpolation on $\O \times [0,T]$
\begin{equation}
\U(x,t) = \frac{t_h - t }{\tauhm} \Hi + \frac{t-t_{h-1}}{\tauhm} \hi \quad  \text{if} \quad x \in \OXi, \quad t\in ]t_{h-1},t_h],
\label{eq:extensionrand}
\end{equation}
and a piecewise approximation 
\begin{equation}
\M (x,t) =\sum_{i=1}^{n} \sum_{h=1}^N \Hi \chi_{]t_{h-1},t_h]} (t) \chi_{\OXi}(x).
\label{eq:piecewiserand}
\end{equation}
Then, we have 
\begin{equation}\tag{$\mathcal{P}^r_n$}
\begin{cases}
\frac{\partial}{\partial t} \check{u}_n(x,t)=-\boldsymbol{\Delta}^{\rand}_p(\bar{u}_n(x,t)), \quad x \in \O, t > 0, \\
\U(x,0) = g_n(x), \quad x \in \O
\end{cases}
\label{neumann11}
\end{equation}
where
\[
g_n(x) = g_i \quad \text{if} \quad x \in\OXi, i \in [n] ,
\]
and the random variable $\rand$ is such that
\[
\rand(x,y) =\Upsilon_{ij} \quad \text{for} \quad (x,y) \in \OXij .
\]


If conditioned with respect to a realization $\bx = (\bx_1,\cdots , \bx_n) $ of the random vector $\bX$, problem~\eqref{neumannrandomdiscrete} can be rewritten on $G_{q_n} (n,K)$ in the following form
\begin{equation}\tag{$\mathcal{P}^{d}_n$}
\begin{cases}
\displaystyle{\frac{\hi - \Hi}{\tauhm}}= \frac{1}{n}\sum\limits_{j = 1}^{n} \lam \abs{u_{j}^{h-1}- \Hi}^{p-2} (u_{j}^{h-1}- \Hi), \quad (i,h) \in [n] \times [N] ,
\\ 
u_{i}^0 = g_i , \quad i \in [n].
\end{cases}
\label{neumanndiscrete}
\end{equation} 

By capitalizing on the results obtained for the the case where $\acc{G_{q_n} (n,K)}_{n \in \N}$ was generated by the deterministic sequence $\bx$, we get the following result.
\begin{theo}
Suppose that $p \in ]1 , + \infty [$, $K \in L^\infty(\O^2)$ is a symmetric and measurable mapping, and $g \in L^{\infty}(\O)$. Let $u$ and $U_h$ denote the unique solutions to~\eqref{neumann} and~\eqref{neumannrandomdiscrete}, respectively. Let $\U$ be the continuous extension of $U_h$ given in~\eqref{eq:extensionrand}. Then, the following hold:
\begin{enumerate}[label=(\roman*)]
\item For $T> 0$, there exists a positive constant $C$, such that for any $\beta > 0$
{\small\begin{equation}
\anorm{u - \U}_{C(0,T; L^p(\O))} \le C\left ( \beta \frac{\log (n)}{n} + \frac{\max\pa{q_n^{-(p-1)}, q_n^{-p/2}}}{  n^{p/2}} \right)^{1/p} + \anorm{K - \Wedg}_{L^p(\O^2)} +  \anorm{g - g_n}_{L^p(\O)} + O(\tau),
\label{eq:estimaterandom}
\end{equation}}
with probability at least $1 - n^{-C \min\{q_n^{2p-1}, q_n^p\}\beta}$. 
\item Suppose furthermore that $g \in \Lip(s, L^q(\O))$ and $K \in \Lip(s', L^q(\O^2))$, $s, s' \in ]0,1]$, and $q_n\norm{K}_{L^\infty(\O^2)} \leq 1$. Let $\theta \eqdef \min\pa{s, s'} \min\pa{1,q/p}$.  Then, for $T > 0$, there exists a positive constant $C$, such that for any $\beta > 0$ and $t \in ]0,e[$
{\small\begin{equation}
\anorm{u - \U}_{C(0,T; L^p(\O))} \le C \pa{\pa{\beta \frac{\log (n)}{n} + \frac{\max\pa{q_n^{-(p-1)}, q_n^{-p/2}}}{ n^{p/2}}^{1/p}} +  \pa{\frac{t\log (n)}{n}}^{\theta}} +O(\tau),
\label{eq:estimateproba2}
\end{equation}}
with probability at least $1 - \bpa{n^{-C \min\{q_n^{2p-1}, q_n^p\}\beta} + 2n^{-t}}$.
\end{enumerate}
\label{theo : totalrandom}
\end{theo}
The dependence of the constant $C$ in the parameters is similar to Remark~\ref{rem:Cparam}(ii).\\

As a preparatory step to prove Theorem~\ref{theo : totalrandom}, the following lemma is instrumental. It establishes that the spacings between the $n$ uniformly distributed nodes are  $O(\log(n)/n)$ with high probability.
\begin{lem}
Consider the sequence of random spacings $(\bX_{(1)}, \bX_{(2)}-\bX_{(1)}, \cdots, 1 - \bX_{(n)})$, where we recall $\acc{\bX_{(i)}}_{i = 1}^{n}$ are the order statistics of $\bX$. Let $t \in ]0,e[$. Then, for any $i \in [n]$
\begin{equation}
\di \eqdef  \bX_{(i)}-\bX_{(i-1)} \leq  t \frac{\log (n)}{n},
\label{eq:diconcentration}
\end{equation}
with probability at least $1 - n^{-t}$.
\label{lem : diconvergence}
\end{lem}
\bpf{of Lemma~\ref{lem : diconvergence}} 
Since $\bX_i $ are i.i.d.~uniform random variables on $\O$, we have, by virtue of~\cite[Theorem~1.6.7]{Reis1989} that the random variables $\di$, $i \in [n]$, have the same distribution as the random variables $Z_i / \sum_{k = 1}^{n+1} Z_k$, where $Z_1, \cdots, Z_{n+1}$ are i.i.d standard exponential random variables.
In addition, invoking \cite[Lemma~1.6.6]{Reis1989}, we know that $S_{n+1}\eqdef \sum_{k = 1}^{n+1} Z_k$ is a Gamma random variable with parameters $(1,n+1)$ (thus having the density $f_{S_{n+1}} (s )= e^{-s} s^n/ n!$, $ s\geq 0$).

Now, combining these two observations, we obtain by straightforward integral calculations that for any $\epsi \in [0,1[$ 
\begin{equation}
\begin{aligned}
\P(\di \geq \epsi )= \P(Z_i \geq \epsi S_{n+1}) &= \P((1- \epsi) Z_i \geq \epsi (S_{n+1} - Z_i))\\
&= \P\left( Z_{n+1}\geq \frac{\epsi}{1- \epsi} S_{n}\right)\\
&= \int_{0}^{+ \infty}  \P\pa{Z_{n+1} \geq  \frac{\epsi}{1- \epsi} s} f_{S_{n}} (s) ds \\
&= \int_{0}^{+ \infty}  e^{- \frac{\epsi}{1- \epsi} s}  e^{-s} \frac{s^{n-1}}{(n-1)!} ds \\
&= (1- \epsi)^n.
\end{aligned}
\end{equation}
The equality of the second line stems from an equality in distribution, since $S_{n+1} - Z_i$ has the same
distribution as $S_n$ and $Z_i$ has the same distribution as $Z_{n+1}$, and the fact that $Z_i$ and $S_{n+1}-Z_i$ are independent.
Taking $\epsi = t \frac{\log (n)}{n} \in ]0,1[$, and using the standard inequality $\log(1-u) \leq -u$, for $u \in [0,1]$, we get 
\[
\P(\di \geq \epsi ) = (1- \epsi)^n = \exp(n \log(1-\epsi)) \leq \exp(-n\epsi) = n^{-t} .
\]
\epf

\bpf{of Theorem~\ref{theo : totalrandom}}
The idea of the proof is to take the conditional probability with respect to a fixed realization $\bx = (\bx_1,\cdots , \bx_n) $ of the random vector $\bX$, then use the bound in Theorem~\ref{theo : convergencedeterministe}, which is independent of $\bx$, and finally integrate with respect to the uniform density on $\O^n$.
\begin{enumerate}[label=(\roman*)]
\item 
We have 
\begin{equation}
\begin{aligned}
\P\pa{ \anorm{u- \check{u}_n}_{C(0,T; L^p(\O))}  \ge \varepsilon'} &= \frac{1}{\abs{\O}^n}\int_{\O^n} \P\pa{\anorm{u- \check{u}_n}_{C(0,T; L^p(\O))}   \ge \varepsilon' | \bX = \bx} d\bx \\
&\leq \frac{1}{\abs{\O}^n}\int_{\O^n} n^{-C  \min\{q_n^{2p-1}, q_n^p\} \beta} d\bx\\
&=n^{-C  \min\{q_n^{2p-1}, q_n^p\}\beta},
\end{aligned}
\end{equation}
with 
\[
\varepsilon' = C\pa{\beta \frac{\log (n)}{n} + \frac{ \max\pa{q_n^{-(p-1)}, q_n^{-p/2}}}{n^{p/2}}}^{1/p} + \anorm{K - \Wedg}_{L^p(\O^2)} +\anorm{g - g_n}_{L^p(\O)} +  O(\tau).
\]
Hence, the desired result, \eqref{eq:estimaterandom} follows from the fact that the obtained estimate in~\eqref{eq:estimate} is uniformally independent of the random choice of $\bx$. 



\item In view of~\eqref{eq: LipK}, we can argue that 
\begin{gather*}
\P\pa{\anorm{K - \Wedg}_{ L^p(\O^2)} \geq \kappa} \leq \P\pa{C(p,q,s')\delta(n)^\theta \geq \kappa} \\
\text{and} \\
\P\pa{\anorm{g - g_n}_{ L^p(\O)} \geq \kappa} \leq \P\pa{C(p,q,s)\delta(n)^\theta \geq \kappa} .
\end{gather*}
Taking $\kappa = \max(C(p,q,s),C(p,q,s'))\pa{t \frac{\log (n)}{n}}^\theta$, for $t \in ]0,e[$, applying Lemma~\ref{lem : diconvergence}, and using a union bound we deduce that the events
\[
\acc{\anorm{K - \Wedg}_{ L^p(\O^2)} \leq \kappa} \qandq \acc{\anorm{g - g_n}_{ L^p(\O)} \leq \kappa}
\]
simultaneously hold with probability at least $1-2n^{-t}$. Denote the events 
\begin{align*}
A_1&: \acc{\anorm{\check{v}_n - \U}_{C(0,T; L^p(\O))}  \leq \epsi} \\
A_2&: \acc{\anorm{K - \Wedg}_{ L^p(\O^2)} \leq \kappa'} \\
A_3&: \acc{\anorm{g - g_n}_{ L^p(\O)} \leq \kappa'}
\end{align*}
and their complements $A_i^c$, where $\epsi = C\pa{\beta \frac{\log (n)}{n} + \frac{ \max\pa{q_n^{-(p-1)}, q_n^{-p/2}}}{n^{p/2}}}^{1/p} + O(\tau)$ and $\kappa' = C \pa{t \frac{\log (n)}{n}}^\theta$, with $C$ the largest constants among the one in claim (i) and $\linebreak\max(C(p,q,s),C(p,q,s'))$. Using again a union bound, we get
\begin{align*}
\P\pa{\anorm{u- \check{u}_n}_{C(0,T; L^p(\O))} \leq \epsi + \kappa'} 
&\geq \P\pa{\cap_{i=1}^3 A_i} = 1 - \P\pa{\cup_{i=1}^3 A_i^c} \\
&\geq 1 - \sum_{i=1}^3 \P\pa{A_i^c} \geq 1 - \pa{n^{-C \min\{q_n^{2p-1}, q_n^p\}\beta} + 2n^{-t}} ,
\end{align*}
which leads to the desired claim.
\end{enumerate}
\epf

\subsection{Rate regimes}
\label{subsec-regimes}
A close inspection of the error bound in~\eqref{eq:estimateproba2} (Theorem~\ref{theo : totalrandom}) reveals three contributions:
\begin{itemize}
\item Spatial discretization: the first contribution is materialized in the first term which scales as (see Remark~\ref{rem:Cparam}(i))
\[
O\pa{\pa{\frac{\log (n)}{n}}^{1/p} + \frac{\max\pa{q_n^{-(1-1/p)}, q_n^{-1/2}}}{n^{1/2}}} .
\]
This term represents the spatial discretization error when approximating the continuous evolution equation~\eqref{neumann} on the random inhomogeneous graph model $G_{q_n}(n,K)$ generated according to~Definition~\ref{def : randomgraph} with the graphon $K$.

\item Data approximation: the second term is $O\pa{\pa{\frac{\log (n)}{n}}^{\theta}}$ which captures the error of discretizting the initial data $g$ and the graphon $K$. The presence of the error on $K$ is clearly tied to the nonlocal nature of the evolution equation on graphs. This approximation error depends on the regularity of $g$ and $K$, and the latter encodes the geometry/structure of the underlying graphs. The more regular $g$ and $K$ are, the faster the convergence rate. 

\item Time discretization: the last term, which is $O(\tau)$, is classical and corresponds to the time discretization error.
\end{itemize} 

At this stage, one may wonder which of the first two terms dominate, or in other words, what are the different regimes exhibited by the convergence rate as a function of the problem parameters $(p,q,s,s')$. This is quite important as it will reveal which nonlocal $p$-Laplacian evolution problems are harder/easier to discretize by highlighting the role of each parameter, and for instance that of $p$ and the impact of nonlocality (i.e.~graphon structure). 

Toward this goal, we first make the error measure in~\eqref{eq:estimateproba2} independent of $p$ and we choose to quantify the error in the classical $L^2(\O)$ norm. Consequently, thanks to Lemma~\ref{lem:quasidist} and Lemma~\ref{inclusion}, as well as boundedness of the solutions, it is not difficult to see that 
{\small\begin{equation}
\anorm{u - \U}_{C(0,T; L^2(\O))} = 
\begin{cases}
O\pa{\pa{\beta \frac{\log (n)}{n}}^{1/p} + \frac{\max\pa{q_n^{-(1-1/p)}, q_n^{-1/2}}}{n^{1/2}} +  \pa{\frac{t\log (n)}{n}}^{\theta} + \tau}, & p \in [2,+\infty[ \\
O\pa{\pa{\beta \frac{\log (n)}{n}}^{1/2} + \frac{\max\pa{q_n^{-(p/2-1/2)}, q_n^{-p/4}}}{n^{p/4}} +  \pa{\frac{t\log (n)}{n}}^{p\theta/2} + \tau^{p/2}} & p \in ]1,2] ,
\end{cases}
\label{eq:estimateprobaL2}
\end{equation}}
holds with probability at least $1 - \bpa{n^{-C \min\{q_n^{2p-1}, q_n^p\}\beta} + 2n^{-t}}$.

To make the rest of the discussion more concrete and also guarantee the convergence of the sequence $\acc{G_{q_n}(n,K)}_{n \in \N}$ to the graphon $K$, we will work under the assumptions of the example in Section~\ref{subsec:graphexample}, i.e. $q_n  = n^{-g(n)}$ with $g(n) \leq c/\log(n)$ for some $c > 0$. Observe that $q_n \in ]0, 1]$, and since $p > 1$, we have
\[
\max\pa{q_n^{-(1-1/p)}, q_n^{-1/2}} \leq q_n^{-1} \leq e^c.
\]
Thus, the second term in~\eqref{eq:estimateprobaL2} reads
\begin{equation}
O\pa{n^{-\min(p/4,1/2)}} .
\label{eq:estimate3}
\end{equation}

Without loss of generality\footnote{This setting is true for many graphons, see, e.g.,~Remark~\ref{rem:Cparam}(iii).}, we also suppose that $s = s'$ and $q \leq p$ so that $\theta = sq/p \in ]0,q/p] \subset ]0,1]$. In this setting, \eqref{eq:estimateprobaL2} reads
{\small\begin{equation*}
\anorm{u - \U}_{C(0,T; L^2(\O))} = O\pa{\pa{\frac{\log (n)}{n}}^{\min(1/p,1/2,sq/p)\min(p/2,1)} + \tau^{\min(p/2,1)}} .
\end{equation*}} 

\begin{figure}[h!]
\centering
\includegraphics[width=0.5\linewidth]{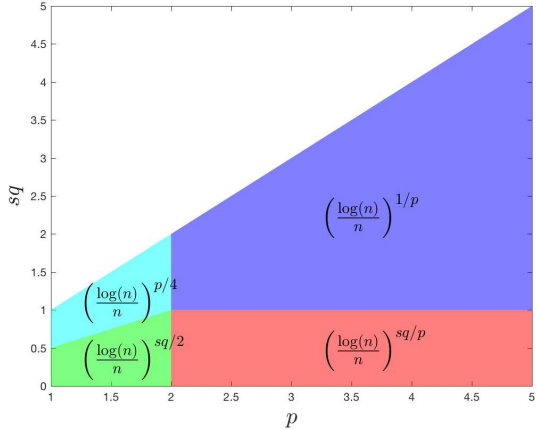}
\caption{Different regimes according to the values of $p$ and $s$, and $q$.}
\label{fig:cell1}
\end{figure}

The term depending on $n$ then exhibits four different regimes as a function of $p$, $s$ and $q$ (see Figure~\ref{fig:cell1}). Indeed, it is straightforward to see that it scales as
\begin{equation*}
\begin{cases}
\pa{\frac{\log(n)}{n}}^{sq/p} \quad \text{for} \quad p \geq 2, \quad sq \in ]0,1],\\ 
\pa{\frac{\log(n)}{n}}^{1/p} \quad \text{for} \quad p \geq 2, \quad sq \in ]1,p],\\ 
\pa{\frac{\log(n)}{n}}^{sq/2} \quad \text{for} \quad p \in ]1,2], \quad sq \in ]0,p/2],\\ 
\pa{\frac{\log(n)}{n}}^{p/4} \quad \text{for} \quad p \in ]1,2], \quad sq \in [p/2,p] .
\end{cases}
\end{equation*}

In particular, the convergence rate shows a transition phenomenon at $p=2$. The rate increases with $p$ for $p \in ]2,+\infty[$ while it decreases with $p$ for $p \in ]1,2]$ and $sq \in [p/2,p]$. As expected, the dependence of the rate on the initial data $g$ and graphon $K$ is more prominent as they become irregular, i.e.~for smaller values of $sq$. For small $sq$ and $p \in ]1,2]$, the rate is independent of $p$.

\appendix

\section{Appendix}
\label{sec:appendix}

\subsection{A key deviation result}
\label{subsec:devZ}

The following lemma establishes a key deviation inequality for $\sup\limits_{t \in [0,T]} \anorm{Z_n(t)}_{p,n}$ where $Z_n(\cdot)$ is the random process defined in~\eqref{zed}.

\begin{lem} Let $Z_n(\cdot)$ be the random process defined in~\eqref{zed}. Then, we have 
\begin{enumerate}[label=(\roman*)]
\item For $p \in ]1, + \infty[$, $T> 0$, there exists a positive constant $C$, such that for any $\beta > 0$
\[
\P \pa{\sup\limits_{t \in [0,T]} \anorm{Z_n(t)}_{p,n} \ge \varepsilon } \le n^{-C \min\pa{q_n^{(2p-1)}, q_n^{p}}  \beta},
\]
with
\[
\varepsilon =  \pa{\beta \frac{\log(n)}{n}+C_3 \max\pa{q_n^{-(p-1)}, q_n^{-p/2}}\frac{1}{n^{p/2}}}^{1/p},
\]
where $C_3$ is a precise constant which will be explicited in the proof.
\item \text{For}  $p \in [2, + \infty[$, suppose that there exists a positive constant $C$, such that for $T > 0$
\[
\displaystyle{\inf\limits_{j \in [n]}}\frac{1}{n} \sum\limits_{i>j} \frac{\alpha^2_{ij}(t)}{q_n}  \wedg (1- q_n \wedg) \ge C.
\]
Then, 
\[
\EE \pa{\anorm{Z_n(t)}_{p,n}^p } \sim \frac{1}{n^{p/2}} .
\]
\end{enumerate}
\label{lem:process}
\end{lem}

To prove this lemma, we need the following deviation inequalities that we include for the reader convenience.

\paragraph{Rosenthal's inequality~\cite{ibragimov1998}.}
Let $n$ be a positive integer, $\gamma \ge 2$ and $U_1,\ldots,U_n$ be $n$ zero mean independent random variables such that $\sup\limits_{i\in [n]}\EE(\abs{U_i}^\gamma )  <\infty$. Then there exists a positive constant $C$ such that
\[
\mathbb{E}\pa{\aabs{\sum_{i=1}^n U_i}^\gamma} \le C \max\pa{ \sum_{i=1}^n\mathbb{E}(|U_i|^\gamma),  \pa{\sum_{i=1}^n\mathbb{E}(U_i^2)}^{\gamma /2}} .
\]

\paragraph{Bernstein's inequality~\cite{massart2007concentration}.}
Let $n$  be a positive integer and $U_1,\ldots,U_n$ be $n$ zero mean independent random variables such that  there exists a positive constant $M$ satisfying $\sup\limits_{i\in [n]}|U_i| \le M <\infty$. Then, for any $\upsilon>0$,
\[
\mathbb{P}\left(\sum_{i=1}^n U_i   \ge \upsilon \right) \le \exp\left(-\frac{\upsilon^2}{2\left(\sum\limits_{i=1}^n\mathbb{E}\left(U_i^2\right)+{\upsilon M}/{3}\right)} \right).
\]

\bpf{of Lemma~\ref{lem:process}} 
\begin{enumerate}[label=(\roman*)]
\item Let us recall that $q_n \lam$ are i.i.d random variables following the Bernoulli distribution with parameter $q_n \wedg$. For the sake of simplicity, set, for $(i,j) \in [n]^2$, $ Y_{ni}\eqdef \abs{\frac{1}{n} \sum\limits_{j=1}^{n} \alpha_{ij} (\lam - \wedg)}^{p}$. We have

\[
I \eqdef \P\pa{ \anorm{Z_n(t)}_{p,n} \ge \varepsilon } 
= \P \pa{  \frac{1}{n} \pa{  \sum\limits_{i=1}^{n}  Y_{ni} - \EE(Y_{ni}) } \ge \varepsilon^p - \frac{1}{n}  \sum\limits_{i=1}^{n} \EE(Y_{ni})}.
\]

It remains to bound $ \EE \pa{ Y_{ni} }$. We distinguish the case when $p \ge 2$ and $p \in ]1, 2[$.
\begin{itemize}
\item $p \ge 2$. Using the Rosenthal inequality with the independent according to $j$ centered random variables $U_{nij} \eqdef \alpha_{ij} (\lam - \wedg)$, we have 
\begin{equation}
\EE \pa{ Y_{ni}}=\frac{1}{n^p} \EE \pa{ \left|\sum_{j=1}^n U_{nij}\right|^p}\le C_1 \frac{1}{n^p}\max\pa{\sum_{j=1}^n\EE(\abs{U_{nij}}^p),  \pa{\sum_{j=1}^n \EE(U_{nij}^2)}^{p /2}}.
\label{rosenthal}
\end{equation}
We have 
\begin{equation*}
\begin{aligned}
\EE \pa{\abs{U_{nij}}^p }&=\frac{\abs{\alpha_{ij}}^p}{q_n^p} (q_n\wedg) (1-q_n\wedg)^p+   \frac{\abs{\alpha_{ij}}^p}{q_n^p} (q_n \wedg)^p (1 -  q_n\wedg) \\
&=\frac{\alpha_{ij}^p}{q_n^{p-1}}   \wedg (1-q_n \wedg)((q_n \wedg)^{p-1}+(1-q_n \wedg)^{p-1}).
\end{aligned}
\end{equation*}
Taking $p = 2$, we get 
\[
\EE (U_{nij}^2 ) =\frac{\alpha^2_{ij}}{q_n}  \wedg (1 - q_n \wedg).
\]
Since $\alpha_{ij} $ and $q_n \wedg$ are both bounded and $p$ being greater than $2$, there exists $C_2>0$, such that, 
\[
\max( \EE\pa{ \abs{U_{nij}}^p}, \EE (U_{nij}^2)) \le C_2  \max \acc{q_n^{-(p-1)}, 1/q_n}.
\] 
Therefore
\begin{equation}
\frac{1}{n} \sum\limits_{i=1}^{n} \EE \pa{ Y_{ni}} \le C_3 \max\pa{q_n^{-(p-1)}, q_n^{-p/2}} \frac{1}{n^{p/2}}. 
\end{equation}
\item $p \in ]1,2[$. With the same steps as above, since $p \in [1, 2[$, applying the Jensen inequality first for the concave function $x \mapsto x^{p/2}$ and second for the convex function $x \mapsto x^2$, we have 
\begin{equation}
\begin{aligned}
 \EE\pa{ Y_{ni}}=\frac{1}{n^p} \EE \pa{ \left|\sum_{j=1}^{n}U_{nij}\right|^p}
& \le \frac{1}{n^p} \pa{\EE\pa{ \left|\sum_{j=1}^{n} U_{nij}\right|^2}}^{p/2}\\
 &\leq  \frac{1}{n^p} \pa{ \sum_{j=1}^{n} \EE \pa{U_{nij}^2 }}^{p/2}\\
 &=\frac{1}{n^p} \pa{\sum_{j=1}^{n}\frac{\alpha_{ij}^2}{q_n}  \wedg (1 - q_n \wedg) }^{p/2}.
\end{aligned}
\end{equation}
Therefore, we have again 
\[
\EE \pa{Y_{ni}}  \le  \frac{C_3}{q_n^{p/2}} \frac{1}{n^{p/2}}.
\]
\end{itemize}
Thus, for any $p > 1$, we get 
\begin{equation}
\frac{1}{n} \sum\limits_{i= 1}^{n} \EE \pa{Y_{ni}} \le C_3 \max\pa{q_n^{-(p-1)}, q_n^{-p/2}}  \frac{1}{n^{p/2}}.
\label{ineq:majoration}
\end{equation}
Hence, setting $W_{ni} = Y_{ni}- \EE\pa{  Y_{ni}} $ and $\lambda =  \varepsilon^p - C_3 \max\pa{q_n^{-(p-1)}, q_n^{-p/2}} \frac{1}{n^{p/2}} $, we have 
\[ 
I \le  \P \pa{ \frac{1}{n} \sum\limits_{i=1}^{n} W_{ni}\ge \lambda }. 
\]
Let $\varepsilon>0$ such that $\lambda > 0$. Observe that the random variables $\{ W_{ni}\}_{i=1}^{n}$ are independent, centred, and obey:
\begin{itemize}[label = $\triangleright $]
\item $\sup\limits_{i \in [n]} \abs{W_{ni}} \le 2\sup\limits_{i \in [n]}  \abs{Y_{ni}}  \le C_{4}$, since $\alpha_{ij}$ and $q_n \wedg$ are both bounded.
\item  $\sum\limits_{i=1}^n\EE\pa{W_{ni}^2}=\sum\limits_{i=1}^n\Var\pa{Y_{ni}}\le\sum\limits_{i=1}^n\EE\pa{Y_{ni}^2}$. Replacing the exponent "$p$" in inequality~\eqref{rosenthal}, by "$2p$" which is greater than $2$, we obtain
\[
\sum\limits_{i=1}^n\EE \pa{
Y_{ni}^2}\le   C_5\max\pa{q_n^{-(2p-1)}, q_n^{-p}} \frac{1}{n^{p-1}} \Rightarrow \sum\limits_{i=1}^n\EE\pa{W_{ni}^2} \le C_5\max\pa{q_n^{-(2p-1)}, q_n^{-p}}\frac{1}{n^{p-1}},
\]
\end{itemize}

We are then in position to apply the Bernstein inequality to $\{ W_{ni}\}_{i=1}^{n}$ according to the index $i$, whence we get, after some elementary algebra
\begin{equation*}
\begin{aligned}
\P \pa{ \frac{1}{n} \sum\limits_{i=1}^{n} W_{ni}\ge \lambda } &\le \exp \pa{ - \frac{n^2 \lambda^2}{ 2  \pa{\sum\limits_{i=1}^{n} \EE \pa{W_{ni}^2}  + n \lambda C_{4} / 3 }}}  \\
&\le \exp \pa{ - \frac{C_6}{2} \min\pa{q_n^{(2p-1)}, q_n^{p}} \frac{n \lambda^2}{n^{-p} + \lambda}}. 
\end{aligned}
\end{equation*}
Taking $\lambda = \beta \frac{\log(n)}{n} > n^{-p}$, for $p > 1$, we have after straightforward calculations 
 \[
 \P \pa{ \frac{1}{n} \sum\limits_{i = 1}^{n} W_{ni}\ge \lambda }  \le  \exp \pa{ - \frac{C_6}{4}  \min\pa{q_n^{(2p-1)}, q_n^{p}}  n \lambda } =  n^{-\frac{C_6}{4}  \min\pa{q_n^{(2p-1)}, q_n^{p}}  \beta}.
 \]
 Therefrom 
 \[
 I \le  \P \pa{  \frac{1}{n} \sum\limits_{i = 1}^{n} W_{ni}\ge  \lambda }  \le  n^{- C  \min\pa{q_n^{(2p-1)}, q_n^{p}} \beta}.
 \]
  For this choice of $\lambda$, observe that 
\begin{equation*}
\begin{aligned}
\lambda=\beta \frac{\log(n)}{n} &\Leftrightarrow   \varepsilon^p-C_3 \max\pa{q_n^{-(p-1)}, q_n^{-p/2}}  \frac{1}{n^{p/2}}=\beta\frac{\log(n)}{n} \\
&\Leftrightarrow  \varepsilon=\pa{\beta \frac{\log(n)}{n}+C_3 \max\pa{q_n^{-(p-1)}, q_n^{-p/2}}\frac{1}{n^{p/2}}}^{1/p}.
\end{aligned}
\end{equation*}
Thus 
\begin{equation}
\P \pa{  \sup\limits_{t \in [0,T]} \anorm{Z_n (t) }_{p,n} \ge \varepsilon } \le  n^{-C \min\pa{q_n^{(2p-1)}, q_n^{p}}\beta}.
\label{ineq:probaZ}
\end{equation}

\item  Set, for $(i,j) \in [n]^2$, $A_n \eqdef \frac{1}{n} \sum\limits_{i=1}^{n} Y_{ni}\eqdef \frac{1}{n} \sum\limits_{i=1}^{n} \abs{Z_{ni}}^p = \frac{1}{n} \sum\limits_{i=1}^{n}\abs{\frac{1}{n} \sum\limits_{j=1}^{n} \alpha_{ij} (\lam- \wedg)}^{p} $.

For $p \in [2, + \infty[$, applying the Jensen inequality twice, we have 
\begin{equation*}
\begin{aligned}
\EE (A_n) = \frac{1}{n} \sum\limits_{i=1}^{n} \EE (Y_{ni}) &= \frac{1}{n^{p+1}} \sum\limits_{i=1}^{n} \EE \pa{ \abs{\sum\limits_{j=1}^{n} U_{nij}}^p} \\
& \ge \frac{1}{n^{p+1}}  \sum\limits_{i=1}^{n}  \pa{\EE \pa{\sum\limits_{j=1}^{n} U_{nij}}^2}^{p/2}\\
&= \frac{1}{n^{p+1}}  \sum\limits_{i=1}^{n}  \pa{\Var \pa{  \sum\limits_{j=1}^{n} U_{nij}  } }^{p/2}\\
&= \frac{1}{n^{p+1}}  \sum\limits_{i=1}^{n} \pa{   \sum\limits_{j=1}^{n} \Var (U_{nij})  }^{p/2}\\
&=  \frac{1}{n^{p+1}}  \sum\limits_{i=1}^{n} \pa{ \sum\limits_{j=1}^{n}  \frac{\alpha^2_{ij}}{q_n}  \wedg (1- q_n \wedg)}^{p/2}\\
&\ge C^{p/2} \frac{1}{n^{p+1}} n n^{p/2}  \ge C_2 \frac{1}{n^{p/2}}.
\end{aligned}
\end{equation*}
Using the mutual independence of  the random variables $\{ \lam \}$ for all $(i,j) \in [n]^2$, 
\begin{equation*}
\begin{aligned}
\EE \pa{   \pa{ \sum\limits_{j=1}^{n} U_{nij}}^2 } &= \Var \pa{    \sum\limits_{j=1}^{n} U_{nij} } \\
\end{aligned}
\end{equation*}
Finally, combined with~\eqref{ineq:majoration}, we conclude that 
\[
C_2 \frac{1}{n^{p/2}} \le \EE \pa{ \anorm{Z_n(t)}_{p,n}^p} \le C_1  \frac{1}{n^{p/2}} .
\]
\end{enumerate}
\epf

\subsection{Approximation theoretic results}
\label{subsec:spaces}
In an effort to make this paper more self-contained we briefly recall some results on functional spaces and approximation theory that our work relies on. But before this, we state the following classical lemma which is useful throughout the paper.

\begin{lem}
For $\alpha \in ]0,1]$ and $a,b \geq 0$, we have
\[
(a + b)^\alpha \leq a^\alpha + b^\alpha .
\]
\label{lem:quasidist}
\end{lem}

\paragraph{$L^p$ spaces embeddings.}
Since $|\O|=1$, we have the classical inclusion $L^q(\Omega) \subset L^p(\Omega)$ for $1 \leq p \leq q < +\infty$. More precisely
\begin{equation}
\anorm{F}_{L^p(\O)} \leq \anorm{F}_{L^q(\O)} \leq \anorm{F}_{L^\infty(\O)} .
\label{classicinclusion}
\end{equation}
We also have the following useful (reverse) bound whose proof is based on H\"older inequality.
\begin{lem}
For any $1 \leq q < p < +\infty$ we have
\[
\anorm{F}_{L^p(\O)} \leq \anorm{F}_{L^\infty(\O)}^{1-q/p}\anorm{F}_{L^q(\O)}^{q/p} .
\]
\label{inclusion}
\end{lem}

\paragraph{Lipschitz spaces $\Lip(s, L^q(\O^d))$~\cite[Ch. 2, \S 6 and 9]{devorelorentz93}.}We introduce the Lipschitz spaces $\Lip(s,L^q(\O^d))$, for $d \in \{1,2\}$, which contain functions with, roughly speaking, $s$ "derivatives" in $L^q(\O^d)$~\cite[Ch.~2, Section~9]{devorelorentz93}. 
\begin{defi}\label{def:lipspaces}
For $F \in L^q(\O^d)$, $q \in [1,+\infty]$, we  define the (first-order) $L^q(\O^d)$ modulus of smoothness by
\begin{equation}
\omega(F,h)_q \eqdef \sup_{\bs z \in \R^d, |\bs z| < h} \pa{\int_{\bs x,\bs x +\bs z \in \O^d}\abs{F(\bs x + \bs z)-F(\bs x)}^q d\bs x}^{1/q} .
\label{modsmooth}
\end{equation}
The Lipschitz spaces $\Lip(s,L^q(\O^d))$ consist of all functions $F$ for which
\[
\abs{F}_{\Lip(s,L^q(\O^d))} \eqdef \sup_{h > 0} h^{-s} \omega(F,h)_q < +\infty .
\]
\end{defi}
We restrict ourselves to values $s \in ]0,1]$ as for $s > 1$, only constant functions are in $\Lip(s,L^q(\O^d))$. It is easy to see that $\abs{F}_{\Lip(s,L^q(\O^d))}$ is a semi-norm. $\Lip(s,L^q(\O^d))$ is endowed with the norm
\[
\anorm{F}_{\Lip(s,L^q(\O^2))} \eqdef  \anorm{F}_{L^q(\O^2)} +  \abs{F}_{\Lip(s,L^q(\O^d))} .
\]
The space $\Lip(s,L^q(\O^2))$ is the Besov space $\mathbf{B}^s_{q,\infty}$~\cite[Ch.~2, Section~10]{devorelorentz93} which are very popular in approximation theory. In particular, $\Lip(1,L^1(\O^d))$ contains the space $\BV(\O^d)$ of functions of bounded variation on $\O^d$, i.e. the set of functions $F \in L^1(\O^d)$ such that their variation is finite:
\begin{equation*}
V_{\O^2}(F) \eqdef \sup_{h > 0}h^{-1}\sum_{i=1}^d\int_{\O^d}\abs{F(\bs x + he_i)-F(\bs x)}d\bs x < + \infty 
\end{equation*}
where $e_i, i \in \{1,d\}$ are the coordinate vectors in $\R^d$; see~\cite[Ch.~2, Lemma~9.2]{devorelorentz93}. Thus Lipschitz spaces are rich enough to contain functions with both discontinuities and fractal structure.

Let us define the piecewise constant approximation of a function $F \in L^q(\O^2)$ (a similar reasoning holds of course on $\O$) on a partition of $\O^2$ into cells $\O_{nij} \eqdef \ens{ ]x_{i-1}, x_i] \times ]y_{j-1}, y_j]} {(i,j) \in [n]^2}$ of maximal mesh size $\delta \eqdef \max\limits_{(i,j) \in [n]^2} \max(\aabs{x_{i} - x_{i-1}},\abs{y_{j} - y_{j-1}})$,
\[
F_n(x,y) \eqdef \sum\limits_{i,j=1}^{n} F_{nij} \chi_{\O_{nij}}(x,y) , \quad F_{ij} = \frac{1}{\abs{\O_{nij}}}\int_{\O_{nij}} F(x,y) dxdy .
\]
Clearly, $F_n$ is nothing but the orthogonal projection of $F$ on the $n^2$-dimensional subspace of $L^{q}(\O^2)$ defined as 
\[
\Span \ens{\chi_{\O_{nij}}}{(i,j)\in [n]^2} .
\]

\begin{lem} 
There exists a positive constant $C_s$, depending only on $s$, such that for all $F \in \Lip(s,L^q(\O^d))$, $d \in \{1,2\}$, $s \in ]0,1]$, $q \in [1,+\infty]$,
\begin{equation}
\anorm{F - F_n}_{L^q(\O^d)} \le C_s \delta^s \abs{F}_{\Lip(s,L^q(\O^d))}.
\label{eq:lipspaceapprox}
\end{equation}
\label{lem:spaceapprox}
\end{lem}

\bpf{} 
Using the general bound \cite[Ch.~7, Theorem~7.3]{devorelorentz93} for the error in spline approximation, and in view of Definition~\ref{def:lipspaces}, we have
\[
\anorm{F - F_n}_{L^q(\O^d)} \leq C_s \omega(F,\delta)_q = C \delta^{s} (\delta^{-s} \omega(F,\delta)_q) \leq C_s \delta^{s} \abs{F}_{\Lip(s,L^q(\O^d))} .
\]
\epf

An immediate consequence is the following result.
\begin{lem}
Assume that $F \in L^{\infty}(\O^d) \cap \Lip (s, L^q(\O^d))$, $d \in \{1,2\}$, $s\in ]0, 1]$, $q \in [1, + \infty]$, and let $p \in ]1, + \infty[$. Then there exists a positive constant $C(p,q,s)$, depending on $p$, $q$ and $s$ such that
\begin{equation}
\anorm{F - F_n}_{L^p(\O^d)} \leq C(p,q,s)\delta^{s \min\{1, q/p\}} .
\label{eq: Lipg}
\end{equation}
\label{lem:estimateglip}
\end{lem}

\bpf{}
We have
\begin{equation*}
\anorm{F - F_n}_{L^p(\O^d)} \leq 
\begin{cases}
\anorm{F - F_n}_{L^q(\O)} \leq C\abs{F}_{\Lip(s,L^q(\O))}\delta^s, \quad \text{if} \quad q\geq p;\\[0.5cm]
\begin{aligned}
\anorm{F - F_n}_{L^{\infty}(\O^d)}^{1- q/p} \anorm{F - F_n}_{L^q(\O^d)}^{q/p} &\leq C\pa{2 \anorm{F}_{L^{\infty}(\O)}}^{1-q/p}\abs{F}_{\Lip(s,L^q(\O^d))}^{q/p} \delta^{s q/p} \\
& \text{otherwise},
\end{aligned}
\end{cases}
\end{equation*}
where we used~\eqref{classicinclusion} (resp.~Lemma~\ref{inclusion}) and Lemma~\ref{lem:spaceapprox} in the first (resp.~second) case.
\epf

\paragraph{Acknowledgement.}
This work was supported by the ANR grant GRAPHSIP. JF was partly supported  by Institut Universitaire de France.

\bibliographystyle{abbrv}
\bibliography{biblioinhomogeneous}

\end{document}